%% file: irreducible_k_omega_arXiv.tex
\magnification=\magstep1 \baselineskip=16pt

\def\Glimm {{\rm Glimm}}

\def\Sub {{\rm Sub}}

\def\Cl{{\rm Cl}}
\input amssym.tex
\input irreduciblequot.ref.tex

\centerline {\bf IRREDUCIBLE QUOTIENT MAPS FROM LOCALLY COMPACT} 
\centerline {\bf SEPARABLE METRIC SPACES}
\medskip
\centerline {\bf A.J. Lazar and D.W.B. Somerset}
\centerline {\bf 17th January 2022}
\bigskip
\bigskip
\noindent {\bf Abstract.} Let $X$ be a Hausdorff quotient of a standard space (that is of a locally compact separable metric space).
It is shown that the following are equivalent: (i) $X$ is the image of an irreducible quotient map from a standard space;
(ii) $X$ has a sequentially dense subset satisfying two technical conditions involving double sequences;
(iii) whenever $q: Y \to X$ is a quotient map from a standard space $Y$, the restriction $q_*|V$ is an 
irreducible quotient map from $V$ onto $X$ (where $q_*:Y_*\to X$ is the pure quotient derived from $q$, and
$V$ is the closure of the set of singleton fibres of $Y_*$). The proof uses extensions of the theorems of Whyburn and Zarikian 
from compact to locally compact standard spaces. The results are new even for quotients of locally compact subsets of the real line.

\bigskip
\noindent Let $X$ and $Y$ be topological spaces and $q:Y\to X$
a surjective map. Recall that $q$ is {\sl irreducible} if $Y$ has no proper closed subset $V$ such that $q(V)=X$, and
$q$ is {\sl inductively irreducible} if $Y$ has a closed subset $W$ such that $q(W)=X$ and $q|W$ is 
irreducible. Irreducible maps have generally been studied under the further hypothesis that $q$ is a closed map [\the\G], [\the\Stricklen]. However, in recent work on
C$^*$-algebras (summarised at the end of the paper) the authors encountered Hausdorff spaces $X$ which were the
images of locally compact separable metric spaces $Y$ under non-closed quotient maps $q$. In some cases $q$ was irreducible,
but in other cases it could be shown that $q$ could not be chosen irreducible. The examples were not esoteric: typically
$Y$ was a locally compact subset of the real line or the plane.

Hausdorff quotients of locally compact separable metric spaces were studied from the 1950s to the 1970s -- with some famous names 
such as Milnor, Steenrod, and Morita -- but they have received little attention since. In an earlier paper, the authors showed
that if $Y$ is a locally compact separable metric space (hereafter a {\sl standard} space) and $q:Y\to X$ is a quotient map
with $X$ Hausdorff then there is an auxiliary standard space $Y_*$ and quotient map $q_*:Y_*\to X$,
and generally $q_*$ has much better properties than $q$. 
Here $Y_*$ is the closure of the set of fibres of $Y$ in the Fell topology on $\Cl(Y)$, the set of closed subsets of $Y$, and for
$F\in Y_*$, $q_*(F)=q(y)$ $(y\in F)$. Then $q_*$ is a {\sl pure} quotient map (Definition 1.2) and this 
permits an extension of Morita's theorem from closed
quotient maps to non-closed quotient maps [\the\BCM; Theorem 4.5]. 

In this paper we characterise those Hausdorff spaces $X$ which are the image of an
irreducible quotient map from a standard space (let us say that such a space is {\sl of irreducible type}). The characterisation is that $X$ should possess a sequentially dense 
subset satisfying two technical conditions involving double sequences (Theorem 5.4). From this it follows (Corollary 5.5) that if $X$ is of irreducible type
then whenever $q: Y \to X$ is a quotient map from a standard space $Y$, there is a unique subset $V\subseteq Y_*$ (namely
the closure of the set of singleton fibres in $Y_*$) such that
the restriction $q_*|V$ is an irreducible quotient map from $V$ onto $X$. The proof involves
extending Whyburn's classic theorem from 1939 (Theorem 2.3) and Zarikian's recent generalisation (Theorem 3.6)
from compact to locally compact standard spaces.
\bigskip
\bigskip
\noindent {\bf 1. Irreducible closed quotient maps}
\bigskip
\noindent We begin by looking at closed quotient maps, which are much simpler than general quotient maps.
We show by a standard argument that every closed quotient map is inductively irreducible (Proposition 1.1)
and we characterise the spaces which are the images of irreducible closed quotient maps from locally compact $\sigma$-compact spaces (Theorem 1.3) and
from standard spaces (Theorem 1.5).
\bigskip
\noindent {\bf Proposition 1.1.} {\sl Let $q : Y\to X$ be a  closed quotient map from a locally compact, $\sigma$-compact Hausdorff space $Y$ to a Hausdorff space $X$.
Then there is a closed subset $Z$ of $Y$ such that $q|Z$ is an irreducible closed quotient map onto $X$.}
 \bigskip
\noindent {\bf Proof.} Let $X_0$ $(X_1)$ be the set of points $x\in X$ such that $q^{-1}(x)$ is non-compact (the boundary of $q^{-1}(x)$ is non-compact). 
Let $X_2$ be the set of isolated points in $X_0$. Then $X_2$ is open and closed [\the\M; Theorem 4] and $X_2\cap X_1 = \emptyset$.
 Also $q^{-1}(X_2)$ is open and closed. Let $Y_0$ be a set consisting of exactly one point from each fibre in $q^{-1}(X_2)$. 
 Then $Y_0$ is closed [\the\M]. Set $Y_1 := q^{-1}(X\setminus X_2)$, a closed and open subset of $Y$. Let $\cal{G}$ be the family of all 
 closed subsets of $Y$ of the form $F\cup Y_0$ where $F$ is a closed subset of $Y_1$ that intersects each fibre of $Y_1$. 
 Order this family by inverse inclusion. A maximal element in this order is the set $Z$ we need. The intersection of a chain in $\cal{G}$ is a closed 
 set of the form $F\cup Y_0$. Then $F$ intersects each fibre $q^{-1}(x)$ with $x\in X\setminus X_0$ by compactness. 
 Thus $q(F\cup Y_0)$ is closed in $X$ and contains a dense subset of $X$; hence $F\cup Y_0$ is in $\cal{G}$. 
 Thus a minimal set $Z$ exists, and $q(Z) = X$, and $q|Z$ is irreducible (minimality) and closed. It is a quotient map since $q$ is closed. Q.E.D.
 \bigskip
\noindent Lasnev [\the\Las] obtained the same conclusion as in Proposition 1.1 under the hypotheses that $Y$ is paracompact and that $X$ is a Fr\'echet-Urysohn space (defined below). 
 
 \bigskip
\noindent {\bf Definition 1.2.} Let $Y$ be a locally compact $\sigma$-compact space and $q:Y\to X$ a quotient map with $X$ Hausdorff.
Recall that $q$ is {\sl pure} [\the\BCM] if there is a subset $D$ of $Y$ 
such that (i) $D$ is dense in $Y$ and the restriction of $q$ to $D$ is injective, and (ii) for every net $(d_{\alpha})$ in $D$, if $d_{\alpha}\to\infty$ (i.e.
eventually escapes from every compact set in $Y$) then $e_{\alpha}\to \infty$ for any other net $(e_{\alpha})$ for
which $q(d_{\alpha})=q(e_{\alpha})$ for all $\alpha$. 

Pure quotient maps have various favourable properties not possessed
by general quotient  maps, and in particular the quotient map $q_*:Y_*\to X$ is always pure. 

 \bigskip
 \noindent Recall that a point $x$ in a topological space $X$ is a $k$-point if $E\subseteq X$ and 
 $x\in \overline E$ implies that there is a compact set $K$ such that $x \in \overline{E \cap K}$.

\bigskip
\noindent {\bf Theorem 1.3.} {\sl Let $q : Y\to X$ be a  pure quotient map from a locally compact, $\sigma$-compact space Hausdorff $Y$ space to a Hausdorff space. 
Let $N$ be the set of all non $k$-points of $X$. The following are equivalent:
 
(1) $N$ is empty;

(2)     there is a closed quotient map from a locally compact, $\sigma$-compact Hausdorff space onto $X$;

(3)      there is an irreducible closed quotient map from a locally compact, $\sigma$-compact Hausdorff space  onto $X$.}
 \medskip
\noindent {\bf Proof.} The equivalence of (1) and (2) follows from [\the\BCM; Theorem 4.1]; (2) implies (3) is Proposition 1.1 above; and 
(3) implies (2) is trivial. Q.E.D.
 
\bigskip
\noindent {\bf Example 1.4.} {\sl The Arens space $S_2$.} This is one of the standard spaces in the subject and will be needed in the
next theorem. Let $Y$ be the disjoint union of a countably infinite collection of convergent sequences.
Let $\{y_n: n\ge 1\}$ be the set of limit points of these sequences and $\{y'_n: n\ge 1\}$ the set of isolated points of the first sequence. 
The non-trivial equivalence classes of $Y$ are the pairs $\{ y'_n, y_{n+1}\}$ $(n\ge 1)$. Let $q:Y\to X$ be the quotient map.
Then $X$ is the Arens space $S_2$ (and we shall call $q:Y\to S_2$ the {\sl standard presentation} for $S_2$). It is easily checked 
that $Y_*$ is homeomorphic to $Y$ and that $q$ is not irreducible. The point $q(y_1)$ is the {\sl base-point} of $S_2$.

\bigskip
\noindent A topological space $X$ is {\sl Fr\'echet-Urysohn} at $x \in X$ if $A\subseteq X$ and $x \in \overline A$ (the closure of $A$) implies the existence of a sequence 
$(x_n)_{n\ge1} \subseteq A$ with $\lim_n x_n = x$. 
If $X$ is Fr\'echet-Urysohn at each point then $X$ is a {\sl Fr\'echet-Urysohn space}. Note that the base-point in $S_2$ is not Fr\'echet-Urysohn.

\bigskip
\noindent {\bf Theorem 1.5.} {\sl Let $q : Y\to X$ be a  pure quotient map from a standard space to a Hausdorff space. 
Let $N$ be the set of all non $k$-points of $X$.
The following are equivalent:
 
(1) $N$ is empty;

(2)      $X$ is a Fr\'echet-Urysohn space;

(3)      there is a closed quotient map from a standard space onto $X$;

(4)      there is an irreducible closed quotient map from a standard space onto $X$;

(5)    $X$ has no closed subset homeomorphic to $S_2$.}
 \medskip
\noindent {\bf Proof.} The equivalence of (1) and (2) follows from  [\the\BCM; Prop. 5.1]. The Arens space $S_2$ contains
a point $x$ which lies in the closure of a subset $C$ from which there is no subsequence converging to $x$. Thus (2) implies (5). The
converse, (5) implies (2), follows from [\the\FST; Proposition 1]. Q.E.D.
\bigskip
\noindent {\bf Examples 1.6.} We close with two examples in which $q$ is irreducible but non-closed. In the first $N$ is a singleton;
in the second, the whole space.
First, let $Y=[0,1)$ and let $(y_n)$ and $(z_n)$ be disjoint sequences of distinct points in $Y$ with $y_n\to 0$
and $z_n\to 1$. The non-trivial fibres are the pairs $(y_n, z_n)$ $(n\ge 1)$. Let $q: Y\to X$ be the
quotient map. Then $q$ is irreducible and $N=\{q(0)\}$.

For the second example, let $Y=(0,1)$ and let $(y_n)$ be a sequence in $(0,1)$ converging to $1$ and
$\{y'_n\}$ a countable dense subset of $Y\setminus \{y_n: n\ge 1\}$. The non-trivial equivalence classes are the pairs 
$\{y_n, y'_n\}$ $(n\ge 1)$. Let $q:Y\to X$ be the quotient map. Then $q$ is irreducible but no point of $X$ is a $k$-point so $N=X$. 

\bigskip
\bigskip
\noindent {\bf 2. Extension of Whyburn's theorem}
\bigskip
\noindent Whyburn's classic theorem states that if $q$ is a continuous surjective map from a compact standard space $Y$ to a Hausdorff
space $X$ then $q$ is irreducible if and only if the set of singleton fibres is dense in $Y$ [\the\Why; Theorem 2].
In this section we extend this result to locally compact standard spaces.
\bigskip
\noindent Let $Y$ be a locally compact, separable metric space and $q:Y\to X$ a continuous surjective map with $X$ Hausdorff.
For $f\in C_0(Y)^+$, define $f^{\sharp}:Y\to {\bf R}$ by $$f^{\sharp}(y)=\sup\{f(w): w\in Y, \ q(w)=q(y)\}.$$ Then $f^{\sharp}$ is upper semi-continuous
on $Y$ and for $f, g\in C_0(Y)^+$ with $\Vert f-g\Vert<\epsilon$, $\sup\{ |f^{\sharp}(y)-g^{\sharp}(y)|: y\in Y\}\le\epsilon$. Thus
if $E$ is a countable dense subset of $C_0(Y)^+$ then $\{f^{\sharp}: f\in E\}$ is a countable dense subset of $\{f^{\sharp}: f\in C_0(Y)^+\}$ with
the supremum norm. It follows that $Y$, being a Baire space, has a dense $G_{\delta}$ subset $G$ consisting of points of continuity for the functions $f^{\sharp}$.
\bigskip
\noindent {\bf Lemma 2.1.} {\sl Let $Y$ be a locally compact, separable metric space and $q:Y\to X$ a continuous surjective map with $X$ Hausdorff. 
Let $y_0\in Y$ and set $F_0=q^{-1}(q(y_0))$. Then $y_0\in G$ (the set of points of continuity for the functions $f^{\sharp}$) if and only if 
whenever $(F_n)_{n\ge 1}$ is a sequence of fibres with limit $F\in \Cl(Y)$ and $y_0\in F$ then $F=F_0$.}
\bigskip
\noindent {\bf Proof.} Suppose first that whenever $(F_n)_{n\ge 1}$ is a sequence of fibres with limit $F\in \Cl(Y)$ and $y_0\in F$ then $F_n\to F_0$.
Let $f\in C_0(Y)^+$ and suppose that $f^{\sharp}(y_0)=\alpha>0$. 
Let $y_m\to y_0$; we may suppose that $(y_m)$ is contained in a neighbourhood of $y_0$ with compact closure. 
There is a subsequence $(y_n)$ such that $(q^{-1}(q(y_n)))$ converges to some $F\in \Cl(Y)$. 
Since $y_0\in F$ we have $F = F_0$ by our assumption. Given any $\epsilon$ with $\alpha > \epsilon > 0$, 
there is $w\in F_0$ such that $f(w) > \alpha - \epsilon$. Passing to a further subsequence if necessary there are 
$w_n\in q^{-1}(q(y_n))$ with $w_n\to w$.  Thus $\lim\sup f^{\sharp}(y_n)\geq \alpha - \epsilon$ and $f^{\sharp}$ is lower semi-continuous at $y_0$.

Conversely, let $(F_n)_{n\ge 1}$ be a sequence of fibres with limit $F\in \Cl(Y)$ and $y_0\in F$ and suppose that there exists $w\in F_0\setminus F$.
Let $f\in C_0(Y)^+$ with $f(w)=1$ and $f$ vanishing in a neighbourhood of $F$. Given $\epsilon>0$, let $K$ be the compact
set $K=\{y\in Y: f(Y)\ge \epsilon\}$. Then eventually $F_n$ misses $K$ so $f^{\sharp}(y)\le\epsilon$ for
$y\in F_n$. But $f^{\sharp}(y_0)\ge 1$. Since there is a sequence
$(y_n)$ with $y_n\in F_n$ and $y_n\to y_0$ we see that $y_0\notin G$. Q.E.D.

\bigskip
\noindent {\bf Lemma 2.2.} {\sl Let $Y$ be a locally compact, separable metric space and $q:Y\to X$ a continuous surjective map with $X$ Hausdorff. 
If $q$ is irreducible, then for every open $U\subseteq Y$ and $n\in {\bf N}$, there is a compact fibre $F$ contained in $U$ with
the diameter of $F$ less than $1/n$.}
\bigskip
\noindent {\bf Proof.} Let $\rho$ be the metric on $Y$ and let $V\subseteq U$ be compact set with non-empty interior $V_0$. Let $y\in V_0$ and
let $U_0=\{w\in V_0: \rho(y,w)<1/2n\}$. Then $U_0$ is a non-empty open set, so $U_0$ contains a fibre $F$ because 
$q$ is irreducible, and $F$ is compact because it is contained in $V$ and has diameter less than $1/n$ because
it is contained in $U_0$. Q.E.D.
\bigskip

\noindent {\bf Theorem 2.3.} {\sl Let $Y$ be a locally compact, separable metric space and
$q:Y\to X$ a continuous surjective map with $X$ Hausdorff. Then the following are equivalent:

(i) $q$ is irreducible;

(ii) for every open $U\subseteq Y$ and $n\in {\bf N}$, there is a compact fibre $F$ contained in $U$ with
the diameter of $F$ less than $1/n$;

(iii) the set of singleton fibres is dense in $Y$ (and indeed contains a dense $G_{\delta}$ subset).}
\bigskip
\noindent {\bf Proof.} (i) implies (ii). This is Lemma 2.2. (iii) implies (i).
If the set of singleton fibres is dense in $Y$ then no proper closed subset of $Y$ can map surjectively onto
$X$, so $q$ is certainly irreducible. 
(ii) implies (iii). Suppose that (ii) holds, and let $G$ be the dense $G_{\delta}$ of Lemma 2.1.
Let $y_0\in G$ with $F_0=q^{-1}(q(y_0))$. Then (ii) implies that there exists a sequence of compact fibres $(F_n)$
with limit $F=\{y_0\}$. By Lemma 2.1, $F_0=F$, so $\{y_0\}$ is a singleton fibre and
the dense $G_{\delta}$ set $G$ consists of singleton fibres. Q.E.D.
\bigskip
\noindent Now let $q:Y\to X$ be an irreducible quotient map. Then each point $y\in Y$ is a limit of a sequence of singleton fibres $(y_n)$,
and $\{y_n\}\to\{y\}$ in the Fell topology, so the singleton $\{y\}$ is an element of $Y_*$. Thus there is a map
$e:Y\to Y_*$ given by $e(y)=\{y\}$. Let $V\subseteq Y_*$ denote the image of $Y$ under $e$. It is easily seen that $V$ is closed
and that $e$ is a homeomorphism of $Y$ onto $V$. Thus $q_*|V$ is an irreducible quotient map from $V$ to $X$. On the other hand,
$q_*:Y_*\to X$ is almost never irreducible. We collect these facts as follows.
\bigskip
\noindent {\bf Corollary 2.4.} {\sl Let $Y$ be a standard space and
$q:Y\to X$ an irreducible quotient map with $X$ Hausdorff.
Then 

(i) $e:Y\to V\subseteq Y_*$ is a homeomorphism;

(ii) $q_*:Y_*\to X$ is irreducible if and only if $q$ is a homeomorphism.}
\bigskip
\noindent {\bf Proof.} (i) Established above. 

(ii) If $q$ is a homeomorphism, then $q_*$ is also a homeomorphism, and hence
irreducible. Conversely, suppose that $q$ is not a homeomorphism. 
Then $Y$ has a non-singleton fibre, so by (i), $V$ is a proper closed subset of $Y_*$
with $q_*|V$ surjective, so $q_*$ is not irreducible. Q.E.D.

\bigskip
\noindent 
Recall the notion of a pure quotient map (Definition 1.2). It is evident from Theorem 2.3 that if $Y$ is a standard space and
$q:Y\to X$ is an irreducible quotient map with $X$ Hausdorff then the set of singleton fibres in $Y$ satisfies the requirements 
for the subset $D$ and hence we have the following.
\bigskip
\noindent {\bf Corollary 2.5.} {\sl Let $Y$ be a standard space and
$q:Y\to X$ a quotient map with $X$ Hausdorff. If $q$ is irreducible then $q$ is a pure quotient map.}

\bigskip
\bigskip
\noindent {\bf 3. Extension of Zarikian's theorem}
\bigskip
\noindent If $q:Y\to X$ is a continuous surjection of compact Hausdorff spaces, then in general
there are many distinct minimal closed sets $K\subseteq Y$ such that $q(K)=X$ (see Proposition 3.3 below). 
Recently, Zarikian [\the\Zar; Theorem 4.2] obtained an extension of Whyburn's theorem, showing that if $q:Y\to X$
is a continuous surjection of compact metric spaces, then the following are equivalent:

(i) there exists a unique minimal closed subset $K\subseteq Y$ such $q(K)=X$;

(ii) the subset of $X$ consisting of the images of singleton fibres in $Y$ is dense in $X$.
\bigskip
\noindent In this section we extend Zarikian's theorem to the case when $Y$ is locally compact and standard (Theorem 3.6). 
We begin with a small point which is trivial when $Y$ is compact, 
but slightly more subtle in the present context. First we need a standard lemma.
\bigskip
\noindent {\bf Lemma 3.1.} {\sl Let $Y$ be a first countable Hausdorff space and $q:Y\to X$ a quotient map with $X$ Hausdorff.  Let $(x_n)$ be a sequence  in $X$.
The following are equivalent:

(i)  $x_n\to x_0\in X$;

(ii) for every subsequence of $(x_n)$ there exists a sub-subsequence $(x_m)$
and $y_m\in q^{-1}(x_m)$ such that $y_m\to y$ for some $y\in q^{-1}(x_0)$.}
\bigskip
\noindent {\bf Proof.} (i) implies (ii).
Suppose that $x_n\to x_0$. If (ii) fails, then there is a subsequence $(x_p)$ with 
no sub-subsequence $(x_m)$ such that there exist $y_m\in q^{-1}(x_m)$ with $y_m\to y$ for some $y\in q^{-1}(x_0)$.
Hence there is a subsequence $(x_p)$ for which
$\bigcup_p q^{-1}(x_p)$ is closed (since $Y$ is first countable and $T_1$ and $X$ is Hausdorff), and hence $\{x_p: p\ge 1\}$ is a closed subset of $X$, contrary to the 
assumption that $x_n\to x_0$. 

(ii) implies (i). Suppose that for every subsequence of $(x_n)$ there exists a sub-subsequence $(x_m)$
and $y_m\in q^{-1}(x_m)$ such that $y_m\to y$ for some $y\in q^{-1}(x_0)$.
Then each subsequence of $(x_n)$ has a sub-subsequence converging to $x_0$; hence $(x_n)$ converges to $x_0$. Q.E.D.
\bigskip
\noindent {\bf Proposition 3.2.} {\sl Let $Y$ be a standard space and
$q:Y\to X$ a quotient map with $X$ Hausdorff. Let $V$ be the closure of the set
of singleton fibres in $Y$ and set $p=q|V$. Then the following are equivalent:

(i) $p$ maps $V$ surjectively (and hence irreducibly) onto $X$;

(ii) the image $E$ of the set of singleton fibres in $Y$ is sequentially dense in $X$.}
\bigskip
\noindent {\bf Proof.} (i) implies (ii). If $p$ is surjective then certainly $p$ is irreducible and $E$
is sequentially dense in $X$. (ii) implies (i). This follows from Lemma 3.1. Q.E.D. 
\bigskip
\noindent Let $q_i: Y_i\to X_i$ $(i=1, 2)$ be two copies of the standard presentation of $X_i=S_2$
(see Example 1.4). Set $Y=Y_1\cup Y_2$ and let $X$ be the quotient space of $Y$ obtained by 
gluing the two copies of $S_2$ at their base-points. Then the set $E$ (the image of the set of
singleton fibres in $Y$) is dense but not sequentially dense in $X$.
\bigskip
\noindent If $q:Y\to X$ is a surjective map, let us say that a subset $W\subseteq Y$ is {\sl full} if $W$ is closed and $q(W)=X$.
Thus $q$ is inductively irreducible if and only if $Y$ has at least one minimal full subset. The authors do not know
of a quotient map $q:Y\to X$ with $Y$ locally compact and Hausdorff and $X$ Hausdorff which is not inductively 
irreducible. 
\bigskip
\noindent {\bf Proposition 3.3.} {\sl Let $Y$ be a locally compact Hausdorff space and $q:Y\to X$ a quotient map with $X$ Hausdorff. 
Suppose either (a) that every fibre in $Y$ is compact or (b) that the set of compact fibres in $Y$ has dense image in $X$ and $q$ is closed.
Then every full subset of $Y$ contains a minimal full subset.}
\bigskip
\noindent {\bf Proof.} Let $K$ be a full subset of $Y$ and $\cal{G}$ be the family of all closed subsets of $Y$ contained in $K$ that intersect each fibre of $Y$.
 Order this family by inverse inclusion. A maximal element in this order is the set $Z$ we need. The intersection of a chain in $\cal{G}$ is a closed 
 set $F$, and $F$ intersects each compact fibre of $Y$ by compactness. 
 Thus $q(F)$ is surjective in case (a) and is closed and contains a dense subset of $X$ in case (b). Hence $F$ is in $\cal{G}$. 
 Thus a minimal set $Z$ exists. Q.E.D.
 
\bigskip
\noindent {\bf Example 3.4.} Let $Y$ be the disjoint union of a countably infinite collection of convergent sequences. 
Each non-trivial fibre consists of an isolated point and a limit point, chosen so that the following hold: (i) each point
except the limit point in the first sequence belongs to a non-trivial fibre; (ii) no fibre contains an isolated point and the limit point from
the same sequence. Let $q:Y\to X$ be the quotient map. Then $X$ is homeomorphic to $S_{\omega}$ (the Arkhangel'skii-Franklin space [\the\AF])
and the only singleton fibre in $Y$ is that consisting of the limit point of the first sequence. The set of limit points $W$ is a minimal full set which maps irreducibly onto
$X$, but $W$ is not unique because any point in $W$ other than the first limit point could be replaced by the other point in its fibre.
\bigskip
\noindent {\bf Lemma 3.5.} {\sl Let $Y$ be a standard space and $q:Y\to X$ a quotient map with $X$ Hausdorff and
let $V$ be a minimal full subset of $Y$. Suppose that $U$ is a non-empty relatively open subset of $V$ with compact
closure $L$ and that $K$ is a compact subset disjoint from $V$ such that $q(K)= q(L)$. Then
$V$ is not  the unique minimal full subset of $Y$.}
\bigskip
\noindent {\bf Proof.} The set $W=K\cup V\setminus U$ is a full subset of $Y$. Let $\cal{G}$ be the family of all full subsets of $Y$ contained in 
$W$, ordered by inverse inclusion. We seek a maximal element in this order. 
The intersection of a chain in $\cal{G}$ is a closed set $F$. Let $v$ be a singleton fibre of $V$ in the complement of $L$.
Then $v$ is a singleton fibre in $W$, so $v\in F$. Since such singleton fibres are dense in $V\setminus L$ (by Theorem 2.3) and $F$ is closed, $F\supseteq V\setminus L$.
Now let $v$ be a singleton fibre of $V$ in $U$. Then $q^{-1}(q(v))\cap W$ is compact, so $q^{-1}(q(v))\cap F$ is non-empty. Let
$v_0\in L$ and let $(v_n)$ be a sequence in $U$ of singleton fibres of $V$ with $v_n\to v_0$ (such a sequence exists
because the singleton fibres of $V$ are dense in $U$, hence in $L$). Let $(w_n)$ be a corresponding sequence in $F$ such
that $q(w_n)=q(v_n)$ for all $n$. Then by the compactness of $F\cap K$, there is a convergent subsequence $(w_m)$ with
limit $w_0\in F\cap K$. Hence $q(w_0)=q(v_0)$, so $q^{-1}(q(v_0))\cap F$ is non-empty. Thus $q(F)\supseteq q(L)$,
so $F$ is a full subset of $Y$. Hence $F$ is in $\cal{G}$, so $\cal{G}$ has a maximal element $Z$, and $Z\not\supseteq U$
so $Z\ne V$. Q.E.D.
 
 \bigskip
 \noindent We now introduce some notation that we need for the next theorem. Let $q:Y\to X$ be a quotient map with $Y$ standard and $X$ Hausdorff.
 Let $Y=\bigcup_n Y_n$ be a compact decomposition of $Y$ and $£V$ a non-empty closed subset of $Y$. For each fixed $n$,
 define a function $d_n: V\to {\bf R}$ as follows.
Let $\rho$ be a metric on $Y$ and for $v\in V$, set $d_n(v)=\sup\{ \rho (y, V): y\in Y_n, q(y)=q(v)\}$ (with $d_n(v)=0$
if $q^{-1}(q(v))\cap Y_n=\emptyset$).  Note that
$d_n(v)$ is finite because $Y_n$ is compact. Clearly $d_n(v)=0$ for all $n$ if and only the fibre $q^{-1}(q(v))$ is contained in $V$.

Let $(v_i)$ be a convergent sequence in $V$ with limit $v_0\in V$, and $(y_i)$ a sequence in $Y_n$ with 
$q(y_i)=q(v_i)$ for all $i$. Then by the compactness of $Y_n$ $(y_i)$ has a subsequence $(y_j)$ convergent to some $y_0$
with $q(y_0)=q(v_0)$, and $\rho (y_0, V)=\lim_j\rho(y_j, V)$. 
Hence the function $d_n$ is upper semi-continuous on $V$. It follows that the
sets $Z_{n,i}=\{ v\in V: d_n(v)\ge 1/i\}$ are closed in $V$.
\bigskip
\noindent The next theorem is an extension of Zarikian's theorem [\the\Zar; Theorem 4.2] .

 \bigskip
\noindent {\bf Theorem 3.6.} {\sl Let $Y$ be a standard space and $q:Y\to X$ a quotient map with $X$ Hausdorff and
let $V$ be a minimal full subset of $Y$. Then the following are equivalent:

(i) $V$ is the unique minimal full subset of $Y$;

(ii) for each $n\ge 1$ and $i\ge 1$, $Z_{n,i}$ has empty interior in $V$;

(iii) $V$ is the closure of the set of singleton fibres of $Y$.}

\bigskip
\noindent {\bf Proof.} (iii) implies (i) is elementary. (i) implies (ii). Suppose that (ii) fails, so that $Z_{{n_0}, {i_0}}$ contains a non-empty relatively open
set $U$ for some $n_0$ and $i_0$. 
By passing to a smaller non-empty relatively open set if necessary, we may assume that $U$ has compact closure $L$
contained in the interior of $Z_{n_0, i_0}$. Then for each $v\in L$, there exists $w\in Y_{n_0}$ with $q(v)=q(w)$ and $\rho(w, V)\ge 1/2i_0$.
Let $Q=\{ w\in Y_{n_0}: \rho(w, V)\ge 1/2i_0\}$. Then $Q$ is a compact subset of $Y$ and $Q\cap q^{-1}(q(L))$ maps surjectively
onto the compact set $q(L)$. Set $p= q|(Q\cap q^{-1}(q(L)))$ and 
let $\alpha: q(L) \to Q\cap q^{-1}(q(L))$ be a cross-section for $p$.
Let $K$ be the closure of $\alpha(q(L))$. Then $K$ is compact, $K\cap L$ is empty, and $q(K)=q(L)$.
It follows from Lemma 3.5 that (i) fails.

(ii) implies (iii).  If each closed subset $Z_{n, i}$ has empty interior in $V$ then $H=V\setminus \bigcup_{n, i} Z_{n, i}$ is a dense $G_{\delta}$ in the 
Baire space $V$, and $H$ consists of the points $v\in V$ whose fibres $q^{-1}(q(v))$ are contained in $V$.
Let $p=q|V$. Then $p$ is irreducible, so the set of points $v\in V$ with singleton fibres $p^{-1}(p(v))$ is also a dense
$G_{\delta}$ subset of $V$ by Theorem 2.3. Intersecting this dense $G_{\delta}$ subset with $H$,
we obtain a dense $G_{\delta}$ subset of $V$ consisting of singleton fibres of $Y$. Q.E.D.
\bigskip
\noindent {\bf Example 3.7.} 
Let $q:Y\to S_2$ be the standard presentation for $S_2$ (see Example 1.4).
The singleton fibres in $Y$ are the limit of the first sequence and the isolated points in all sequences but the first.
The closure $V$ of this set of singleton fibres is the unique minimal full subset and the restriction map $p=q|V$ is irreducible but not
a quotient map.
\bigskip
\bigskip
\noindent {\bf 4. Irreducible quotient maps from subsets of $Y_*$}
\bigskip
\noindent In this section we begin the problem of determining the Hausdorff spaces $X$ which are of irreducible type  (i.e. the image of an 
irreducible quotient map from a standard space).
Following the lead of Corollary 2.4, given a quotient map $q:Y\to X$
with $Y$ standard and $X$ Hausdorff, we focus on the question of whether there is a closed subset $V\subseteq Y_*$
such that $q_*|V$ is an irreducible quotient map onto $X$. We begin by identifying a necessary condition.

Let $X$ be a Hausdorff space which is a quotient of a standard space. We say that $x\in X$ is a point of {\sl compact type} 
if there is a compact subset $K\subseteq X$ such that whenever $(x_n)$ is a sequence in $X$ converging to $x$, eventually $x_n$ belongs to $K$.
The next proposition explains the name.
\bigskip
\noindent {\bf Proposition 4.1.} {\sl Let $q: Y \to X$ be a pure quotient map from a standard
space $Y$ to a Hausdorff space $X$. Let $x\in X$. Then $x$ is of compact type if and only if $q^{-1}(x)$ is compact.}
\bigskip
\noindent {\bf Proof.} Suppose first that $q^{-1}(x)$ is compact. Let $C$ be a compact subset of $Y$ containing 
$q^{-1}(x)$ in its interior and set $K=q(C)$. Let $(x_n)$ be sequence in $X$ converging to $x$. Then by Lemma 3.1, for any subsequence of $(x_n)$ we may find
a sub-subsequence $(x_m)$ and $y_m\in q^{-1}(x_m)$ such that $y_m\to y$ for some $y\in q^{-1}(x)$. Hence eventually
$y_m\in C$ so eventually $x_m\in K$. Thus eventually $x_n\in K$. (Note that we have not used the purity of $q$.)

Conversely, suppose that $q^{-1}(x)$ is non-compact and let $K$ be any compact subset of $X$. Let $Y=\bigcup_{i\ge 1} Y_i$ be a compact
decomposition for $Y$ and let $D$ be the distinguished dense subset of Y (Definition 1.2). Then by standard theory, there exists $i_0\ge 1$ such that 
$K\subseteq q(Y_{i_0})$. By [\the\BCM; Lemma 3.2], there exists $i_1\ge i_0$ such that 
$$D_{i_0}=\{d\in D: q(d)\in q(Y_{i_0})\}\subseteq Y_{i_1}.$$
Let $y\in q^{-1}(x)\setminus Y_{i_1}$ and let $(d_n)$ be a sequence in $D\setminus Y_{i_1}$ converging to $y$. Then $q(d_n)\notin q(Y_{i_0})$
for every $n$, so
$q(d_n) \notin K$, but $q(d_n)\to q(y)=x$. Q.E.D.
\bigskip
\noindent {\bf Corollary 4.2.} {\sl Let $q: Y \to X$ be a quotient map from a standard
space $Y$ to a Hausdorff space $X$. Then the set of points of non-compact type in $X$ is an $F_{\sigma\delta}$.}
\bigskip
\noindent {\bf Proof.} 
Let $Y_*=\bigcup_n Y_n$ be a compact decomposition for $Y_*$ and set $U_n=Y_*\setminus Y_n$ 
and $V_n=\overline{U_n}$. Then by Proposition 4.1, the set of points of non-compact type in $X$ is equal to $\bigcap_n q_*(V_n)$, and each $q_*(V_n)$ is a countable union
of compact sets, being the continuous  image of a closed subset of the $\sigma$-compact space $Y_*$. Q.E.D.
\bigskip

\noindent The authors do not know of an example in which the set of points of non-compact type is not an $F_{\sigma}$.

It follows from Proposition 4.1 (and Corollary 2.5) that if $q: Y \to X$ is an irreducible quotient 
map from a standard space $Y$ to a Hausdorff space $X$ then $x$ is of compact type whenever $q^{-1}(x)$ is a singleton. 
Hence, by Theorem 2.3, a necessary condition for $X$ to be of irreducible type is that the set of 
points of compact type be sequentially dense in $X$. Not all Hausdorff quotients of standard spaces satisfy this, as we now show.
\bigskip
\noindent {\bf Example 4.3.} Let $Y$ be the disjoint union of a countably infinite collection of convergent sequences.
Each non-trivial fibre consists of an isolated point and a
countably infinite number of limit points, chosen so that the following hold: (i) each point 
except the limit point of the first sequence belongs to a non-trivial fibre; (ii) no fibre contains an isolated point and the limit point from
the same sequence. Let $q:Y\to X$ be the quotient map. Taking $D$ to be the set of isolated points of $Y$, we see that $q$ is a pure quotient map, and hence Proposition 4.1 applies.
Every fibre in $Y$ except one is non-compact so only one point of the infinite space $X$ is of compact type. Hence $X$ is not of irreducible type
(although $Y$ has many minimal full subsets, so $X$ is the image of an irreducible {\sl map} from a standard space).
\bigskip
\noindent The next example shows that even when every point of $X$ is of compact type, $X$ may still fail to be of irreducible
type.
\bigskip
\noindent {\bf Example 4.4.} Let $Y=\{0\}\cup \{ 1/n: n\ge 1\}\cup [2,3)$. The non-singleton fibres are the pairs $(1/n, 3-1/n)$ $(n\ge 1)$.
Let $q:Y\to X$ be the quotient map. Then $X$ is Hausdorff and all the points of
$X$ are of compact type by the proof of Proposition 4.1 (first paragraph). Set $x_0=q(0)$ and $x_n=q(1/n)$.
Then $x_n\to x_0$, and the points $x_n$ are non-isolated in $X$, but for any sequences of distinct points
$(x_{i, n})_{i\ge 1}$ with $x_{i,n}\to x_n$, we have that no diagonal sequence from the family
$\{ x_{i, n}: i, n\ge 1\}$ converges to $x_0$. If $X$ were the image of a standard space $Y_1$ under an irreducible quotient
map, the dense subset of singleton fibres in $Y_1$ would imply that such convergent
diagonal sequences did exist. Thus $X$ is not of irreducible type.

\bigskip
\noindent The third example shows that when $q$ is irreducible, $X$ may have a dense subset of points
of non-compact type.
\bigskip
\noindent {\bf Example 4.5.} Let $Y=(0,1]$ and let $(y_n)$ be a sequence of distinct points in $Y$ 
converging to $0$. Let $\{z_n\}$ be a countable dense subset of $Y\setminus \{y_n: n\ge 1\}$.
The non-trivial equivalence classes each consist of one point $z_n$ and
infinitely many points from the sequence $(y_n)$, with every $z_n$ in such an equivalence class. 
Then $q:Y\to X$ is irreducible (hence pure, by Corollary 2.5) and $q(z_n)$ is not of compact type for each $n$, so the set of points
of non-compact type is dense in $X$.
\bigskip
\noindent 
For a standard space $Y$, we now define an order on $Y_*$ by saying that
for $z_1, z_2\in Y_*$, $z_1\le z_2$ if $z_2$ is a subset of $z_1$ when both are viewed as
closed subsets of $Y$. In this order, the minimal elements of $Y_*$ are the fibres of the 
map $q:Y\to X$, and the set of minimal elements is a dense subset of $Y_*$ by definition. 
More generally, if $V$ is a closed subset of $Y_*$ then a standard argument shows that
each $v\in V$ dominates one or more minimal elements of $V$.

The existence of maximal elements in  $Y_*$ is less certain. If $x\in X$ is a point of compact type
then the fibre $q_*^{-1}(x)$ in $Y_*$ is compact, and hence each element in $q_*^{-1}(x)$ is dominated by a maximal
element. For points $x$ of non-compact type, however, the presence of maximal elements 
in $q_*^{-1}(x)$ depends on $Y$, as the following example shows.
\bigskip
\noindent {\bf Example 4.6.} Let $Y=\bigcup_{n\ge 1} Y_n$ be the disjoint union of subspaces $Y_n$, each
consisting of $n$ distinct sequences (the first sequence, the second sequence, etc.)
converging to the same limit point $y_0^n$. The equivalence classes, which are all non-trivial, are
(1) the set $\{y_0^n: n\ge 1\}$ of limit points; (2) for an isolated point $y\in Y$-- say $y$ is the $i$th term in the $j$th sequence
in $Y_k$ -- the equivalence class consists of all other points which are the $i$th term in the $j$th sequence in some $Y_n$ $(n\ge j)$. 
Then the nested sets $L_k=\{ y_0^n: n\ge k\}$ $(k\ge 1)$ are limits of sequences of fibres in $\Cl(Y)$ and $Y_*$ is homeomorphic to 
the Arens space $S_2$ with the base-point deleted. Setting $x_0=q(y_0^1)$, the only non-trivial fibre in $Y_*$ is $q_*^{-1}(x_0)$ 
which consists of the set of limit points and has no maximal elements. Note that $q_*$ is irreducible, and that in $Y_{**}$, $q_{**}^{-1}(x_0)$
does have maximal elements.

\bigskip
\noindent In the next example, $Y_*$ has no maximal elements at all.
\bigskip
\noindent {\bf Example 4.7.} Let $Y_1=[0,1)$ and $Y_{n+1}=Y_n\times Y_1$ $(n\ge 1)$. Let $Y$ be the disjoint union of the
spaces $Y_n$. Define an equivalence relation on $Y$ by identifying $y\in Y_n$ with $(y,0)\in Y_{n+1}$ at each stage. Let $q:Y\to X$
be the quotient map. Then $Y_*$ has no maximal elements.
Taking $D=Y_1\cup (Y_2\setminus Y_1\times \{0\}) \cup\ldots\cup (Y_{n+1}\setminus Y_n\times \{0\}) \cup\ldots$, we
see that $q$ is pure and hence that $X$ has no elements of compact type by Proposition 4.1.
\bigskip
\noindent {\bf Proposition 4.8.} {\sl Let $q: Y \to X$ be a quotient map from a standard
space $Y$ to a Hausdorff space $X$. Let $V$ be a full subset of $Y_*$ and suppose
that $q_*|V$ is a quotient map. Then 

(i) for every $z\in Y_*$ there exists $v\in V$ such that $z\le v$;

(ii) $V$ contains the maximal elements of $Y_*$, and each maximal element of $V$ is a maximal element of $Y_*$.}
\bigskip
\noindent {\bf Proof.} Let $(f_n)$ be a sequence of minimal elements of $Y_*$ converging to $z$.
Then $x_n=q_*(f_n)$ converges to $x=q_*(z)$, so the closed sets $q_*^{-1}(x_n)\cap V$
cannot converge to infinity in the Fell topology on $\Cl(Y_*)$ (Lemma 3.1). Hence there is a sequence $(y_n)$ 
with $y_n\in q_*^{-1}(x_n)\cap V$ for each $n$, such that $(y_n)$ has a convergent
subsequence $(y_m)$ with limit $v$, say. Then $y_m\ge f_m$ so $v\ge z$.  Thus (i) holds, and
(ii) is immediate from (i). Q.E.D.

\bigskip
\noindent Proposition 4.8 shows that the closure $V$ of the set of maximal elements of $Y_*$
is important in the study of irreducible quotient maps, and we now consider $p=q_*|V$. We have already noted that 
each point of compact type in $X$ belongs to $p(V)$.
\bigskip
\noindent {\bf Lemma 4.9.} {\sl Let $q: Y \to X$ be a quotient map from a standard
space $Y$ to a Hausdorff space $X$ and suppose that the set $E$ of points of compact type is sequentially dense in $X$. 
Let $V$ be the closure in $Y_*$ of the set of maximal elements of $Y_*$ and set $p=q_*|V$. 
Suppose that whenever $(x_n)$ is a sequence in $E$ converging to a limit in $X$, 
the sequence of closed sets $p^{-1}(x_n)\not\to\infty$ in the Fell topology. Then $p$ is surjective.}
\bigskip
\noindent {\bf Proof.}  Let $(x_n)$ be a convergent sequence in $E$ with limit $x\in X$. By assumption there
is a sequence $(v_n)$, with $v_n\in p^{-1}(x_n)$ for each $n$, such that $(v_n)$ has a convergent subsequence with limit $v$.
Then $p(v)=x$, and thus $p$ is surjective. Q.E.D.

\bigskip
\noindent {\bf Theorem 4.10.} {\sl Let $q: Y \to X$ be a quotient map from a standard
space $Y$ to a Hausdorff space $X$ and suppose that the set $E$ of elements of compact type is sequentially dense in $X$. 
Let $V$ be the closure in $Y_*$ of the set of maximal elements of $Y_*$ and set $p=q_*|V$. Then $p$ is a quotient map onto $X$ if and only if whenever
$(x_n)$ is a convergent sequence in $X$, the sequence of closed sets $p^{-1}(x_n)\not\to\infty$ in the Fell topology.}
\bigskip
\noindent {\bf Proof.} Suppose first that $p$ is a quotient map onto $X$ and that $(x_n)$ is a convergent sequence with limit $x_0\in X$. If the sequence of closed 
sets $p^{-1}(x_n)\to\infty$ in the Fell topology, then $(x_n)$ is non-constant and $x_0\ne x_n$ for all $n$. The
set $Z=\bigcup_n p^{-1}(x_n)$ is a closed saturated subset of $V$, so $p(Z)= \{x_n: n\ge 1\}$ is closed in $X$ since $p$ is
a quotient map. This contradiction shows that the condition is necessary. 

Conversely, suppose that the condition holds (Lemma 4.9 shows that $p$ is surjective,
so the condition makes sense).
We show that $p$ is a quotient map. Let $S$ be a subset of $X$ with $p^{-1}(S)$ closed in $V$. Let $(z_m)$ be a sequence in $q_*^{-1}(S)$
converging to a limit $z_0$. We wish to show that $z_0\in q_*^{-1}(S)$.
Set $x_n=q_*(z_n)$ and $x_0=q_*(z_0)$. Then $(x_n)$ converges to $x_0$. 
By the assumed condition, there is a sequence $(v_n)$, with $v_n\in p^{-1}(x_n)$ for each $n$, 
such that $(v_n)$ has a convergent subsequence with limit $v$. Then $v\in p^{-1}(S)$ and $p(v)=x_0$ so $x_0\in S$.
Hence $z_0\in q_*^{-1}(S)$, and thus $S$ is closed and $p$ is a quotient map. Q.E.D.
\bigskip
\noindent Combining Proposition 4.8 and Theorem 4.10, we see that if $X$ and $Y$ satisfy the conditions of Theorem 4.10, then $V$
is the unique minimal closed subset of $Y_*$ such that $p=q_*|V$ is a quotient map. We now show that $p$ has to be irreducible
in this case. We need a preliminary proposition.
\bigskip
\noindent {\bf Proposition 4.11.} {\sl Let $q: Y \to X$ be a quotient map from a standard
space $Y$ to a Hausdorff space $X$ and suppose that the set $E$ of elements of compact type is sequentially dense in $X$. 
Let $V$ be the closure of the set of maximal elements of $Y_*$ and set $p=q_*|V$. Suppose that
$p$ is a quotient map. Then for each maximal $f_0\in V$ there is a sequence of elements $(w_n)$ in $V$
such that $p^{-1}(p(w_n))\to \{ f_0\}$ in the Fell topology on $\Cl(V)$. }
\bigskip
\noindent {\bf Proof.} Let $f_0$ be a maximal element of $V$ (and hence a maximal element of $Y_*$
by Proposition 4.8) corresponding to a non-empty closed subset $F_0$ of $Y$. 
Then there is a sequence $(f_n)$ of minimal elements of $Y_*$ converging to $f_0$, corresponding to a sequence $(F_n)$ of fibres
of $Y$ converging to $F_0$ in the Fell topology on $\Cl(Y)$. The maximality of $f_0$ in $Y_*$ is equivalent
to the minimality of the set $F_0$ as a limit of fibres of $Y$. Set $x_n=q_*(f_n)$ $(n\ge 0)$. Then $x_n\to x_0$, so 
no subsequence of $(p^{-1}(x_n))$ converges to $\infty$ in the Fell topology on $\Cl(V)$ by Theorem 4.10. Indeed any subsequence of $(p^{-1}(x_n))$
contains a convergent subsequence, and since $f_0$ is maximal in $Y_*$, and for each $v_n\in p^{-1}(x_n)$, $v_n\ge f_n$, the only possible limit
set is $\{f_0\}$. Hence  $p^{-1}(x_n)\to \{ f_0\}$ in the Fell topology. Q.E.D.
\bigskip
\noindent {\bf Theorem 4.12.} {\sl Let $q: Y \to X$ be a quotient map from a standard
space $Y$ to a Hausdorff space $X$ and suppose that the set of elements of compact type is sequentially dense in $X$. 
Let $V$ be the closure in $Y_*$ of the set of maximal elements of $Y_*$ and set $p=q_*|V$. Suppose that
$p$ is a quotient map. Then $p$ is irreducible.}
\bigskip
\noindent {\bf Proof.}
 Let $E$ be the dense $G_{\delta}$ subset of $V$ 
consisting of the points of continuity arising from the continuous surjective map $p:V\to X$ (Lemma 2.1). Let $v_0$ belong to $E$. We
show that the fibre $F_0=p^{-1}(p(v_0))$ of $v_0$ in $V$ is a singleton, from which the result follows.

Let $(U_n)$ be a decreasing neighbourhhood base for $v_0$ in $V$ and suppose that each $U_n$ has compact closure $K_n$
contained in $U_{n-1}$ (and $U_1$ has compact closure $K_0$). 
Let $(v_n)$ be a sequence of maximal elements of $V$ converging to $v_0$ with $v_n\in U_n$ for each $n$
and let $V=\bigcup_n V_n$ be a compact decomposition for $V$. Then by Proposition 4.11, applied to each $v_n$ in turn, for each $n$ there
exists $w_n\in V$ with fibre $F_n=p^{-1}(p(w_n))$ such that $w_n\in U_n$ and 
$F_n$ is disjoint from the compact set $V_n\setminus U_n$. It follows that $(w_n)$ converges to $v_0$ and that $(F_n)$ converges
in the Fell topology on $\Cl(V)$ to $F=\{v_0\}$. By Lemma 2.1, $F_0=F$. Q.E.D.

\bigskip
\noindent It will follow from work in the next section that $V$ in Theorem 4.12 is the the unique minimal full subset of $Y_*$,
and hence is the closure of the set of singleton fibres of $Y_*$.
\bigskip
\bigskip
\noindent {\bf 5. Images of irreducible quotient maps}
\bigskip
\noindent
In this section we characterise the Hausdorff spaces $X$ which are of irreducible type, showing
that if $X$ is such a space and $q:Y\to X$ is any quotient map from a standard space
to $X$ then $p=q_*|V$ is an irreducible quotient map, where $V$ is the closure of the set of singleton fibres in $Y_*$.

\bigskip
\noindent Suppose then that $q: Y \to X$ is an irreducible quotient map from a standard space $Y$ to a Hausdorff space $X$
and let $E$ be the image of the set of singleton fibres in $Y$. Then $E$ is sequentially dense in $X$ by Theorem 2.3,
and we have seen that $E$ consists of points of compact type. Furthermore, it is not difficult to see (Lemma 5.1) that $E$ and $X$ have the following two properties:
\medskip
\noindent {\bf Property (a).} Whenever $(x_m)$ is a sequence of distinct points in $E$ with limit $x_0\in X$ there is a subsequence
$(x_n)$ such that for all sequences $(x_n^j)_j$ in $X$ with $x_n^j\to x_n$ there exist natural numbers
$k_n$ $(n\ge 1)$ such that the set $\{ x_n^j: n\ge 1, j\ge k_n\}$ has compact closure in $X$.
\medskip
\noindent {\bf Property (b).} Whenever $(x_m)$ is a sequence of distinct points in $X$ with limit $x_0\in X$ there is a subsequence $(x_n)$ and 
sequences $(x_n^j)_j$ in $E$ with $x_n^j\to x_n$ such that every diagonal sequence from the set $\{ x_n^j: n, j\ge 1\}$ converges to $x_0$.
\bigskip
\noindent We shall see that the possession of a sequentially dense subset $E$ of points of compact type with properties (a) and (b)
characterises the spaces $X$ that are of irreducible type.
\bigskip
\noindent {\bf Lemma 5.1.} {\sl Let $q: Y \to X$ be an irreducible quotient map from a standard
space $Y$ to a Hausdorff space $X$ and let  $E$ be the image in $X$ of the set of singleton fibres of $Y$. 
Then $E$ is sequentially dense in $X$ and consists of points of compact type, and $E$ and $X$ have properties (a) and (b).}
\bigskip
\noindent {\bf Proof.} By Theorem 2.3, the set $F$ of singleton fibres is dense in $Y$. Set $E=q(F)$.
Then $E$ is sequentially dense in $X$, and points in $E$ are of compact type by Proposition 4.1.

Let $(x_m)$ be a sequence of distinct points in $E$ with limit $x_0\in X$. Then by Lemma
3.1, there is a subsequence $(x_n)$ and $y_n\in q^{-1}(x_n)$ such that $y_n\to y$ for some $y\in q^{-1}(x_0)$.
Let $K$ be a compact neighbourhood of $y$ containing the sequence $(y_n)$ in its interior. Then since
each set $q^{-1}(x_n)$ is a singleton contained in $K$, Lemma 3.1 implies that any
sequence $(x_n^j)$ in $X$ with $x_n^j\to x_n$ eventually lies in $q(K)$. Hence
there exist natural numbers
$k_n$ $(n\ge 1)$ such that the set $\{ x_n^j: n\ge 1, j\ge k_n\}$ has compact closure in $X$. Thus property (a) holds. 

Let $(x_m)$ be a sequence of distinct points in $X$ with limit $x_0\in X$. Then by Lemma 3.1, there is a convergent
sequence $(y_n)$ in $Y$ with limit $y_0\in q^{-1}(x_0)$ such that $(x_n=q(y_n))$ is a subsequence of $(x_m)$. 
Let $\rho$ be a metric on $Y$. By passing to a subsequence if need be, we may assume that $\rho (y_0, y_n)<1/n$ for each $n\ge 1$.
For each $n\ge 1$, let $(y_n^j)$ be a sequence of singleton fibres 
converging to $y_n$ with $\rho (y_n, y_n^j)<1/(nj)$ for each $j\ge 1$. 
Set $x_n^j=q(y_n^j)$. Then  every diagonal sequence from the set $\{ x_n^j: n, j\ge 1\}$ converges to $x_0$. 
Thus property (b) holds. Q.E.D.

\bigskip
\noindent We note for later (proof of Theorem 6.5) that $E$ in Lemma 5.1 can be replaced by any subset of $E$ that is sequentially dense in $X$. 
Next we interpret property (a) in $Y_*$.

\bigskip
\noindent {\bf Lemma 5.2.} {\sl Let $q: Y \to X$ be a quotient map from a standard
space $Y$ to a Hausdorff space $X$. Let $E$ be a sequentially dense subset of $X$ such that
each $x\in E$ is of compact type. The following are equivalent:

(a1) whenever $(x_m)$ is a sequence of distinct points in $E$ with limit $x_0\in X$ there is a subsequence
$(x_n)$ such that for all sequences $(x_n^j)$ in $X$ with $x_n^j\to x_n$ there exist natural numbers
$k_n$ $(n\ge 1)$ such that the set $\{ x_n^j: n\ge 1, j\ge k_n\}$ has compact closure in $X$;

(a2) whenever $(x_m)$ is a sequence of distinct points in $E$ with limit $x_0\in X$ there is a subsequence
$(x_n)$ and a compact set $K\subseteq Y_*$ such that $q^{-1}_*(x_n)\subseteq K$ for all $n\ge 1$.}
\bigskip
\noindent {\bf Proof.}  Suppose first that (a2) does not hold. Let $Y_*=\bigcup_m Y_m$ be a compact decomposition for $Y_*$. Then 
we may find a convergent sequence $(x_m)$ of distinct points in $E$ with limit $x_0\in X$ such that
$q_*^{-1}(x_{m})\not\subseteq Y_m$.   Let $y_{m}\in 
q_*^{-1}(x_{m})\setminus Y_m$ and let $(f_{m}^j)$ be a sequence of minimal elements of $Y_*$ converging to
$y_{m}$. We may assume that $f_{m}^j\notin Y_m$. Then $y_{m}\to \infty$ and any diagonal sequence $(f_n)_n$ from the
set $\{ f_{m}^j: m, j \ge 1\}$ also converges to infinity. Hence any sequence $(g_n)_n$ with $g_n\in Y_*$ and
$g_n\ge f_n$ also converges to infinity. Set $x_{m}^j=q_*(f_{m}^j)$. Then any diagonal sequence or subsequence from the
set $\{ x_{m}^j: i, j \ge 1\}$ is non-convergent by the above, so there does not exist a subsequence
$(x_n)$ and natural numbers $k_n$ $(n\ge 1)$ such that the set $\{ x_n^j: n\ge 1, j\ge k_n\}$ has compact closure in $X$.
Thus (a2) implies (a1).

Conversely, suppose that (a2) holds. Let $(x_m)$ be a sequence of distinct points in $E$ with limit $x_0\in X$. By assumption, there is a subsequence
$(x_n)$ and a compact set $K\subseteq Y_*$ such that $q^{-1}_*(x_n)\subseteq K$ for all $n\ge 1$. By covering $K$ with a finite number of open sets with compact
closure and taking $K$ to be the union of these compact sets, we may suppose that each $q^{-1}_*(x_n)$ is contained in the interior of $K$. 
Let  $(x_n^j)$ in $X$ with $x_n^j\to x_n$. It follows from Lemma 3.1 that for each $n\ge 1$, there exists $k_n$ such that $q^{-1}_*(x_n^j)$ meets $K$ for all $j\ge k_n$.
For $j\ge k_n$ let $y_n^j\in q^{-1}_*(x_n^j)\cap K$. Then $q_*(K)$ is a compact
subset of $X$ containing the set $\{ x_n^j=q_*(y_n^j): n\ge 1, j\ge k_n\}$, which thus has compact closure. Q.E.D.
\bigskip
\noindent Lemma 5.2 shows that the easiest way to understand property (a) is that its failure implies the existence of a convergent sequence $(x_n)$ in $E$
such that $q_*^{-1}(x_{m})\not\subseteq Y_m$ where $Y_*=\bigcup_m Y_m$ is any given compact decomposition for $Y_*$. Thus property (a) is excluding this possibility.
\bigskip
\noindent We also need the following implication of property (b) in $Y_*$.
\bigskip
\noindent {\bf Lemma 5.3.} {\sl Let $q: Y \to X$ be a quotient map from a standard
space $Y$ to a Hausdorff space $X$. Suppose that $X$ has a sequentially dense subset $E$ with property (b) above.
Set $F=\{q^{-1}(x): x\in E\}$ (regarded as a subset of $Y_*$) and let $w_0$ be a maximal element of $Y_*$. 
Then there is a sequence $(f_n)\subseteq  F$ with $f_n\to w_0$.}
\bigskip
\noindent {\bf Proof.} 
Let $Y=\bigcup_n Y_n$ be a compact decomposition for $Y$, and for each $n$, set $Z_n=\{ z\in Y_*:
 q^{-1}(q_*(z))\cap Y_n\ne\emptyset\}$. Then $Y_*=\bigcup_n Z_n$ is a compact decomposition for $Y_*$
 and if $g\in Z_n$ and $f\le g$ then $f\in Z_n$. For each $n$, set $K_n=q_*(Z_n)$.
 
If $w_0$ is a minimal element of $Y_*$ then $\{w_0\}$ is a singleton fibre in $Y_*$.
In this case, either $w_0$ is isolated in $Y_*$, so that $x_0$ is isolated in $X$ and hence belongs to $E$,
or $w_0$ is non-isolated, so there is a sequence of distinct points (which can be assumed to minimal)
converging to $w_0$. Thus, for general $w_0$,  we may assume that $(g_m)$ is a sequence of distinct minimal elements of $Y_*$ with $g_m\to w_0$. Set
$x_m=q_*(g_m)$ and $x_0=q_*(w_0)$. Then by property (b)
 there is a subsequence $(x_n)$ and sequences $(x_n^j)$ in $E$ with $x_n^j\to x_n$ 
 such that every diagonal sequence from the set $\{ x_n^j: n, j\ge 1\}$ converges to $x_0$.

Let $K = \{x_n^j: n\geq 1, j\geq n\}$, with closure $\overline K$. Then $\overline K=K\cup \{x_n:n\geq 0\}$ and is sequentially compact.
Hence $\overline K\subseteq K_l$ for some $l$, and we pick 
$f_n^j \in K_l$ such that $q_*(f_n^j) = x_n^j$, $f_n^j\in F$, $n\geq 1$, $j\geq n$. By the compactness of $K_l$, we may find a convergent
subsequence $(f_n^i)\to$ with limit $w_n$, $w_n\geq g_n$. 
Every subsequence of $(w_n)$ has a subsequence that converges to $w_0$ by the maximality of $w_0$, and hence $w_n\to w_0$. 
Thus we may find a diagonal sequence of $(f_n^i)$ converging to $w_0$.  Q.E.D.

\bigskip
\noindent We are now ready for the main characterisation.
\bigskip
\noindent {\bf Theorem 5.4.} {\sl Let $q: Y \to X$ be a quotient map from a standard
space $Y$ to a Hausdorff space $X$. Then the following are equivalent:

(i) $X$ has a sequentially dense subset $E$ such that each $x\in E$ is of compact type and 
$X$ and $E$ have properties (a) and (b);

(ii) the map $q_*|V: V\to X$ is an irreducible quotient map from $V$ (the closure of the set of singleton fibres of $Y_*$)
onto $X$.}

\bigskip
\noindent {\bf Proof.} (ii) implies (i) is Lemma 5.1.

(i)$\Rightarrow$(ii).
Let $V$ be the closure of the set of maximal elements of $Y_*$ and set $p=q_*|V$.
First we show that $p$ is surjective. If $x_0\in X$ is isolated then $x_0\in E\subseteq p(V)$, so we may suppose that $x_0$ is non-isolated. Let 
$(x_m)$ be a sequence of distinct points in $E$ converging to $x_0$. Then by property (a) and Lemma 5.2
there is a subsequence $(x_n)$ and a compact set $K\subseteq Y_*$ such that $q^{-1}_*(x_n)\subseteq K$ for all $n\ge 1$. Let
$v_n\in q^{-1}_*(x_n)\cap V \subseteq K$.
Then $(v_n)$ has a convergent subsequence, with limit $v$, say, in $V\cap K$ so $p(v)=x_0$ and $p$ is surjective. 

To show that $p$ is irreducible, let $\rho$ be a metric on $Y_*$ and let $U$ be an open subset of $V$. Let $w_0$ be a maximal element of $Y_*$
belonging to $U$. Set $F=\{ q^{-1}(x): x\in E\}$ (regarded as a subset of $Y_*$),
and let $(f_m)$ be a sequence of elements of $F$ converging to $w_0$ (by Lemma 5.3, and thus using property (b)).
Set $x_m=q_*(f_m)$. Since each $q_*^{-1}(x_m)$ is compact, 
there is a corresponding sequence $(w_m)$ of maximal elements of $W$ with $q_*(w_m)=x_m$ for each $n\ge 1$.
Then $(x_m)$ converges to $x_0=q_*(w_0)$, so by property (a), there is a subsequence
$(x_n)$ and a compact set $K\subseteq Y_*$ such that $q^{-1}_*(x_n)\subseteq K$ for all $n\ge 1$.
Hence the subsequence $(w_n)\subseteq K$, so by passing to a further subsequence if necessary, we may
suppose that $(w_n)$ is convergent, and its limit must be $w_0$ by the maximality of $w_0$ in $Y_*$. Thus
eventually $w_n$ belongs to $U$.

Suppose for a contradiction that the diameters of the fibres $q_*^{-1}(q_*(f_n))$  do not converge to zero.
Then there exists $\epsilon >0$ and a subsequence $(x_i)$ of $(x_n)$ with $v_i\in q_*^{-1}(q_*(f_i))$ and
$\rho(v_i, f_i)\ge\epsilon$ for all $i\ge 0$. By passing to a further subsequence, we may suppose that $(v_i)$
is convergent, and again its limit must be $w_0$ by the maximality of $w_0$ in $Y_*$. This contradicts
the triangle law for metrics. 
Thus the diameters of the fibres $q_*^{-1}(q_*(f_n))$ converge to zero. This implies two things. 

First, the diameters of the fibres $p^{-1}(q_*(f_n))$ also converge to
zero, so eventually these compact fibres are contained in $U$. Thus $p$ is irreducible by Theorem 2.3. 
Secondly, because $(w_n)$ is convergent to $w$, and the diameters of the fibres $q_*^{-1}(q_*(f_n))$
converge to zero, and $U$ and $w$ were arbitrary, it follows that the sets $Z_{n, i}$ of Theorem 3.6 have empty interior. 
Thus, since $p$ is irreducible, $V$ is the unique minimal full subset of $Y_*$, and $V$ is the closure of the set
of singleton fibres of $Y_*$.

Finally, we show that $p$ is a quotient map. Let $Y_*=\bigcup_n Y_n$ be a compact decomposition of $Y_*$
and let $x_n\to x_0$ in $X$. By Theorem 4.8, we want to show that $(p^{-1}(x_n))$ does not converge to infinity. 
Let the sequences $(x_n^j)_j$ be given by property (b), and set $A = \{x_n^j|n\geq 1, j\geq n\}$. 
We claim that  there is $l\ge 1$ such that $A\subseteq p(Y_l\cap V)$. Assume not. 
Then for each $n\ge 1$, there exists $x_n'\in A$ with $x_n'\notin p(Y_n\cap V)$.  Then $(x_n')$ has a converging subsequence, say $(x'_m)$. 
There exists $l\ge 1$ such that $q_*^{-1}(x'_m)\subseteq Y_l$ for all $n$ by property (a2). Thus
$p^{-1}(x'_m)\subseteq Y_l\cap V$ for each $m$, a contradiction. Thus $\{x_n\}\subseteq \overline A\subseteq p(Y_l\cap V)$ for some $l$, so
$(p^{-1}(x_n))$ does not converge to infinity. Q.E.D.
\bigskip
\noindent {\bf Corollary 5.5.} {\sl Let $X$ be a Hausdorff quotient of a standard space. Then the following are equivalent:

(i) $X$ is the image of an irreducible quotient map from a standard space;

(ii) $X$ has a sequentially dense subset $E$ such that each $x\in E$ is of compact type and 
$X$ and $E$ have properties (a) and (b);

(iii) whenever $q: Y \to X$ is a quotient map from a standard space $Y$, the restriction $q_*|V:V\to X$ is an 
irreducible quotient map from $V$ (the closure of the set of singleton fibres of $Y_*$)
onto $X$.}
\bigskip
\noindent As one application, we see that the Arens space $S_2$ is not of irreducible type because if
$q:Y\to S_2$ is the standard presentation (Example 1.4) then $q_*|V$ maps surjectively onto
$S_2$ but is not a quotient map. 
\bigskip
\bigskip
\noindent {\bf 6. Some applications}
\bigskip
\noindent In this section we give some further applications of Theorem 5.4 and Corollary 5.5.
\bigskip
\noindent We begin by picking up a point from Section 4. We saw there that if $q: Y\to X$
is a quotient map from a standard space $Y$ to a Hausdorff space $X$, and if
$X$ has a sequentially dense subset of elements of compact type and
$p:=q_*|W$ is a quotient map onto $X$ (where $W$ is the the closure of the set of maximal elements of $Y_*$) 
then $p$ is irreducible (Theorem 4.12). However, this left open the question
of whether $W$ was the unique minimal full subset of $Y_*$ or equivalently whether $W$
coincided with $V$, the closure of the set of singleton fibres of $Y_*$. With Theorem 5.4, we immediately know the answer:
$p$ is an irreducible quotient map, so $X$ is the image of an irreducible quotient, so
$V$ is the unique minimal full subset in $Y_*$ and hence $V=W$.

Let us say that a subset $Z\subseteq Y_*$ has {\sl the quotient property} if $Z$ is full and $q_*|Z$ is a quotient
map. It follows from Theorem 5.4 that if $X$ is of irreducible type
then $Z\subseteq Y_*$ has the quotient property if and only if $Z$ contains $V$, the closure of the
set of singleton fibres. When $X$ is not of irreducible type the order-structure of the set
of subsets of $Y_*$ with the quotient property can be very different: for instance, if $X=S_2$ (with the standard presentation), one can easily
find a descending chain of subsets of $Y\cong Y_*$ with the quotient property whose intersection does not have that property. Perhaps
this phenomenon occurs more widely but we have not pursued this subject.
\bigskip
\noindent We now give an application of Theorem 5.4 to one of the most tractable classes of Hausdorff quotients of a standard space,
namely those for which the quotient map can be assumed to be closed. We saw in  Proposition 1.1 that if $q:Y\to X$ is a closed quotient map with
$Y$ standard and $X$ Hausdorff then $Y$ has a minimal full subset $Z$ such that $q|Z$ is a quotient map, so $X$ is of irreducible type.
Furthermore, if $q$ is closed then $q_*:Y_*\to X$ is also closed [\the\BCM; Proposition 3.7]. Thus we have the following.
\bigskip
\noindent {\bf Corollary 6.1.} {\sl Let $q:Y\to X$ be a quotient map with $Y$ standard and $X$ Hausdorff. 
Let $V$ be the closure of the set of singleton fibres of $Y_*$ and set $p=q_*|V$.
If $q_*$ is closed (in particular if $q$ is closed) then $V$ is the unique minimal full subset of $Y_*$
and $p$ is an irreducible quotient map.}
\bigskip
\noindent The most important special case of Corollary 6.1 is when $Y$ is a compact standard space, but here a simpler proof is available
(which will also be useful in Proposition 6.4).
\bigskip
\noindent {\bf Theorem 6.2.} {\sl Let $X$ be the Hausdorff quotient of a compact standard space $Y$. 
Let $V$ be the closure of the set of singleton fibres of $Y_*$ and set $p=q_*|V$. Then $V$ is the unique minimal full subset of $Y_*$
and $p$ is an irreducible quotient map.}
\bigskip
\noindent {\bf Proof.} Let $\psi : X\to Y_*$ be defined by $\psi(x) := q^{-1}(x)$ (regarded as a point in $Y_*$). 
Then $\psi$ is continuous at $x\in X$ if and only if $q_*^{-1}(x)$ is a singleton (for otherwise 
there exists $(x_n)$ in $X$  with $x_n\to x$ and $(q^{-1}(x_n))$ converging in the Fell topology to a proper
subset of $q^{-1}(x)$).
 We show that there is a dense $G_{\delta}$ subset $B$ of $X$ such that $\psi$ is continuous at the points of $B$. The result will then
 follow by [\the\Zar; Theorem 4.2]. 
 
First note that if $x\in X$ is isolated then $\psi$ is continuous at $x$.
For $f\in C(Y)^+$ define $f^{\flat}(x) := \sup\{f(y)\mid y\in q^{-1}(x)\}$, $x\in X$; then $f^{\flat}$ is upper semicontinuous on $X$. 
If $\|f_1 - f_2\| < \epsilon$ then $\|f_1^{\flat} - f_2^{\flat}\| < \epsilon$. Let $(f_n)$ be a dense sequence in $C(Y)^+$. 
There is a dense $G_{\delta}$ subset, say $B$, of $X$ at the points of which all of the functions $f_n^{\flat}$, 
and consequently all the functions $f^{\flat}$ for $f\in C(Y)^+$, are continuous.
 
The claim is that $\psi$ is continuous at the points of $B$. Let $x$ be a non-isolated point of $B$ and 
suppose that $(x_m)$ is a sequence in $X$ that converges to $x$ such that $(\psi(x_m))$ converges in $Y_*$ to some 
$E\subsetneq \psi(x)$. Let $y\in \psi(x)\setminus E$ and let $U$ be an open subset of $Y$ that contains $E$ and 
such that $y\notin \overline{U}$. Let $f\in C(Y)^+$ satisfy $0\leq f\leq 1$, $f(y) = 1$, $f|{\overline{U}} = 0$. 
Given $\epsilon > 0$, for $m$ large enough we have $f^{\flat}(x_m) < \epsilon$ (the complement of $U$ is compact) 
while $f^{\flat}(x) = 1$, a contradiction to the continuity of $f^{\flat}$ at $x$. Q.E.D. 
\bigskip

 \noindent Let $X$ be the image of a quotient map $q$ from a standard space $Y$ then $X$ and let $L$ be the
 open subset of $X$ consisting of elements which have a compact neighbourhood. If $q$ is closed then $L$ is a dense 
 subset of $X$ with discrete complement [\the\M; Theorem 4]. Note that in this case, $q^{-1}(L)$ must be dense in $Y$, so $L$ is
 sequentially dense in $X$. When $L$ is dense in $X$ but $X\setminus L$ is non-discrete, $L$ may fail to be sequentially
 dense, as the following example shows.
 \bigskip
\noindent {\bf Example 6.3.} Let $Y=(0,1)\cup (1,2)$ and let $(y_n)$ be a sequence in $(0,1)$ converging to $1$,
and $\{ y'_n\}$ a countable dense subset of $(1,2)$. The non-trivial equivalence classes in $Y$ are the pairs
$\{y_n, y'_n\}$ $(n\ge 1)$. Let $q:Y\to X$ be the quotient map. Then the set of singleton fibres is dense in $X$,
so $q$ is irreducible, and the set $L$ of points of local compactness is dense in $X$ but not sequentially dense.
 \bigskip
\noindent  Our final application is to those spaces $X$ for which $L$ is sequentially dense in $X$. 
 We show that Theorem 5.4 considerably simplifies for these spaces. First we need a general proposition, extending
 Theorem 6.2.
 \bigskip
 \noindent {\bf Proposition 6.4.} {\sl Let $Y$ be a standard space and $q:Y\to X$ a quotient map
 with $X$ Hausdorff. Let $L$ be the set of points in $X$ which have a compact neighbourhood. Then
 the set $E$ of points $x\in L$ for which $q_*^{-1}(x)$ is a singleton is dense in $L$.}
 \bigskip
\noindent  {\bf Proof.} Let $V$ be an open subset of $X$ with compact closure $W\subseteq L$. 
 Let $Y=\bigcup_n {Y_n}$ be a compact decomposition of $Y$. There is $Y_m$ such that $q(Y_m)$ contains $W$, and by renumbering 
 we may assume that $m=1$. The restriction of $q$ to $Y_n\cap q^{-1}(W)$ is $q_n$. Then $q_n:(Y_n\cap q^{-1}(W))\to W$ is a quotient
map.  It follows from the proof of Theorem 6.2 that there is a dense $G_{\delta}$ subset $A$ of $V$ 
such that each function $x\to q_n^{-1}(x)$ is continuous at each $x\in A$. 
 
 Let $x_0\in A$. We claim that the function $x\to q^{-1}(x)$ is continuous at $x_0$. Suppose not. Then there is a sequence $(x_k)$ in $V$ 
 that converges to $x_0$ such that $(q^{-1}(x_k))$ converges to some proper subset $B\subseteq q^{-1}(x_0)$ in the Fell topology. There exists 
 $y \in q^{-1}(x)$ with a compact neighbourhood $C$ disjoint from $B$. Eventually $q^{-1}(x_k)$ is disjoint from $C$. 
 Let $n$ be such that $Y_n\supseteq C$. Then the function $x\to q^{-1}_n(x)$  is not continuous at $x_0$, a contradiction. 
 Hence for $x\in A$, $q_*^{-1}(x)$ is a singleton, so $x\in E$. Q.E.D.
 \bigskip
\noindent  Elementary examples show that $E$ varies with $q$ and that membership of $E$ is not an intrinsic property
of a point of $X$ (by contrast with the property of being a point of compact type).
\bigskip
\noindent Now let $q:Y\to X$ be a quotient map from a standard space $Y$ to a Hausdorff space $X$ and suppose that $L$ is sequentially dense in $X$
(note that this property is independent of $Y$ [\the\BCM; Corollary 4.2]). Let $E$ be the dense subset of $L$ consisting of points $x$ for which $q_*^{-1}(x)$ is a singleton. 
Then $E$ consists of points of compact type, and $E$ automatically has property (a) (using Lemma 3.1 and property (a2) in Lemma 5.2).
However, the example of $X=S_2$ shows that $E$ need not be sequentially dense in $X$ or have property (b).
\bigskip
 \noindent {\bf Theorem 6.5.} {\sl Let $Y$ be a standard space and $q:Y\to X$ a quotient map
 with $X$ Hausdorff. Suppose that $L$ is sequentially dense in $X$. Let $E$ be the set of points in $L$ whose
 inverse image under $q_*$ is a singleton. Then the following are equivalent:
 
 (i) $X$ is the irreducible quotient of a standard space;
 
 (ii) $E$ is sequentially dense in $X$ and has property (b);
 
 (iii) for all $x_0\in X\setminus L$, whenever $(x_m)$ is a sequence of distinct points in $X$ with limit $x_0\in X$ there is a subsequence $(x_n)$ and 
sequences $(x_n^j)$ in $E$ with $x_n^j\to x_n$ such that every diagonal sequence from the set $\{ x_n^j: n, j\ge 1\}$ converges to $x_0$.}
\bigskip
\noindent {\bf Proof.} (i) implies (ii). Suppose that $X$ is the irreducible quotient of a standard space.
Let $F$ be the set of singleton fibres in $Y_*$, $V$ the closure of the set of $F$, $p$ the restriction of
$q_*$ to $V$, $U$ the inverse image of $L$ under $p$, and $F'$ the set of singleton fibres in $U$ (so
that $p(F')=E$. Then $U$ is relatively open because $L$ is open, and $U$ is dense in $V$ because $L$ is
sequentially dense in $X$ and $p$ is a quotient map. Hence $F'$ must be dense in $V$, so $E$
is sequentially dense in $X$. It follows that $E$ has property (b) by the remark after Lemma 5.1.

Conversely, (ii) implies (i) by Theorem 5.4 and the remarks preceding this theorem. Thus (i) and (ii) are equivalent.

For (ii)$\Leftrightarrow$(iii) note that $L$ is metrizable and $E$ is dense in $L$ so $E$ automatically has property (b)
with regard to points $x_0\in L$. Thus (iii) is the part of property (b) that can fail; and if (iii) holds, $E$ is sequentially
dense in $X$. Q.E.D.

\bigskip
\noindent {\bf Example 6.6.} Let $C\subseteq [0,1]$ be the Cantor set and let $Y=C\cup [2,3)$. The non-trivial
equivalence classes in $Y$ are the pairs $\{ c, c+2\} (c\in C\setminus \{1\})$. Let $q:Y\to X$ be the quotient map.
Then $X=L\cup \{q(1)\}$, but the set $E$ of points in $L$ whose inverse image under $q$ is a singleton is not sequentially dense in $X$,
and property (b) fails at $q(1)$. Hence $X$ is not of irreducible type.
\bigskip
\noindent {\bf Remark 6.7.} We conclude by indicating the connection of this paper with C$^*$-algebras. 
Let $A$ be a C$^*$-algebra and let $Id(A)$ and $Id'(A)$ be sets of all closed two-sided ideals of $A$, and
all proper closed two-sided ideals of $A$, respectively. The main topology on $Id(A)$ is the weak
topology generated by the functions $I\to\Vert a+I\Vert$ ($I\in Id(A)$, $a\in A$) which is compact and Hausdorff.
One important subset of $Id'(A)$ is $\Sub(A)$ which is the closure of the set of minimal closed primal ideals of $A$.
Thus $\Sub(A)$ is locally compact, being compact if $A$ is unital, $\sigma$-compact if $A$ is $\sigma$-unital, and a
standard space if $A$ is separable. Furthermore $\Sub(A)$ carries a natural equivalence relation so there is a
quotient map $\psi: \Sub(A)\to \Glimm(A)$, and if $A$ is $\sigma$-compact, the quotient space $\Glimm(A)$ is a 
Tychonov space whose Stone-Cech compactification $\beta\Glimm(A)$ is the character space of the centre of the multiplier algebra of $A$.

The authors' original interest was to determine the possible Glimm spaces for particular classes of separable C$^*$-algebras, such as
postliminal, liminal, sub-homogenous, and others. One common constraint was that $\Sub(A)$ often has a dense set of singleton fibres,
and as the quotient map $\psi$ need not be closed, this led to an unresearched area of topology not far off the beaten track.
In [\the\BCM] we studied the general structure of Hausdorff quotients of locally compact $\sigma$-compact Hausdorff spaces, extending Morita's theorem 
from closed to non-closed quotients; and in this paper we have focussed on standard spaces and irreducible quotient maps. We hope to
give applications to C$^*$-algebras in subsequent work.
\bigskip
\centerline{\bf References}
\bigskip

\input TempReferences
\bigskip
\noindent {\sl School of Mathematical Sciences, Tel Aviv University, Tel Aviv 69978, Israel }

\noindent {\sl e-mail:} aldo@tauex.tau.ac.il, douglassomerset@yahoo.com

\end

%% file: irreduciblequot.ref.tex
\newcount\refno\refno=0
\chardef\other=12
\newwrite\reffile
\immediate\openout\reffile=TempReferences
\outer\def\ref{\par\medbreak\global\advance\refno by 1
 \immediate\write\reffile{}
 \immediate\write\reffile{\noexpand\item{[\the\refno]}}
 \copytoblankline}
\def\copytoblankline{\begingroup\setupcopy\copyref}
\def\setupcopy{\def\do##1{\catcode`##1=\other}\dospecials
 \catcode`\|=\other \obeylines}
{\obeylines \gdef\copyref#1
 {\def\next{#1}%
 \ifx\next\empty\let\next=\endgroup %
 \else\immediate\write\reffile{\next} \let\next=\copyref\fi\next}}
 
\ref A.V. Arhangel'skii and S.P. Franklin, Ordinal invariants for topological spaces, Michigan Math. J., 15 (1968), 313-320.

\newcount\AF\AF=\refno

\ref S.P. Franklin and B.V. Smith Thomas, A survey of $k_{\omega}$-spaces, Topology Proceedings, 2 (1977), 111-124.

\newcount\FST\FST=\refno

\ref G. Gruenhage, Irreducible restrictions of closed mappings, Topology and its Applications, 85 (1998), 127-135.

\newcount\G\G=\refno

\ref N. Lasnev, Continuous decompositions and closed mappings of metric spaces, Soviet Math. Dokl., 6 (1965), 1504-1506.

\newcount\Las\Las=\refno

\ref A.J. Lazar and D.W.B. Somerset, Pure quotients and Morita's Theorem for $k_{\omega}$-spaces, Canadian Mathematical Bulletin 
(to appear).

\newcount\BCM\BCM=\refno

\ref K. Morita, On closed mappings, Proc. Japan Acad. 32 (1956), 538-543.

\newcount\M\M=\refno

\ref S.A. Stricklen, Jr., Closed mappings of nowhere locally compact metric spaces, Proc. Amer. Math. Soc., 68 (1978), no. 3, 369-374.

\newcount\Stricklen\Stricklen=\refno

\ref G.T. Whyburn, On irreducibility of transformations, Amer. J. Math., 61 (1939), 820-822.

\newcount\Why\Why=\refno

\ref V. Zarikian, A generalization of Whyburn's theorem, and aperiodicity for abelian C$^*$-inclusions, https://arxiv.org/abs/2011.13460v1

\newcount\Zar\Zar=\refno

\immediate\closeout\reffile